\documentclass[b5paper,10pt,leqno]{article}
\usepackage{graphicx}
\usepackage{amsmath,amsfonts}
\usepackage{amsthm}

\theoremstyle{plain}
\newtheorem{theorem}{Theorem}[section]
\newtheorem{lemma}[theorem]{Lemma}

\newtheorem{corollary}[theorem]{Corollary}
\newtheorem*{theorem-nn}{Theorem}
\newtheorem*{corollary-nn}{Corollary}

\theoremstyle{definition}

\newtheorem{example}[theorem]{Example}

\newcommand{\bZ}{\mathbb{Z}}

\newcommand{\bQ}{\mathbb{Q}}
\newcommand{\opi}{\overline{\pi}}

\begin{document}
\begin{center}
{\Large\bf On the simplest quartic fields and related Thue equations}
\mbox{}\\[11pt]
{\large\sc
Akinari Hoshi
}
\mbox{}\\[5pt] 
Department of Mathematics, Rikkyo University\\
3--34--1 Nishi-Ikebukuro Toshima-ku Tokyo, 171--8501, Japan\\
{\tt hoshi@rikkyo.ac.jp}\\[5pt] 


\footnotetext{
This work was partially supported by Rikkyo University Special Fund for Research.}
\end{center} 

\begin{abstract}
Let $K$ be a field of char $K\neq 2$. 
For $a\in K$, we give an explicit answer to the field isomorphism problem 
of the simplest quartic polynomial $X^4-aX^3-6X^2+aX+1$ over $K$ 
as the special case of the field intersection problem via multi-resolvent polynomials. 
From this result, over an infinite field $K$, we see that the polynomial gives the same 
splitting field over $K$ for infinitely many values $a$ of $K$. 
We also see by Siegel's theorem for curves of genus zero that 
only finitely many algebraic integers $a\in\mathcal{O}_K$ in a number field $K$ 
may give the same splitting field. 
By applying the result over the field $\mathbb{Q}$ of rational numbers, we establish a 
correspondence between primitive solutions to the parametric family of quartic Thue equations
\[
X^4-mX^3Y-6X^2Y^2+mXY^3+Y^4=c,
\]
where $m\in\mathbb{Z}$ is a rational integer and $c$ is a divisor of $4(m^2+16)$, 
and isomorphism classes of the simplest quartic fields. 
\end{abstract}

\section{Introduction and main results}

Let $K$ be a field of char $K\neq 2$ and $K(s)$ the rational function field over $K$ 
with variable $s$. 
We take the simplest quartic polynomial 
\[
f_s(X):=X^4-sX^3-6X^2+sX+1\in K(s)[X]
\]
with discriminant $4(s^2+16)^3$. 
The Galois group $\mathrm{Gal}_{K(s)} f_s(X)$ of the polynomial $f_s(X)$ over $K(s)$ 
is isomorphic to the cyclic group $C_4$ of order four. 

In the case where $K=\bQ$, for $a\in\bZ\setminus\{0,\pm 3\}$ the polynomials $f_a(X)$ 
are irreducible over $\bQ$ with $\mathrm{Gal}_\bQ f_a(X)\cong C_4$ and the splitting 
fields $\mathrm{Spl}_\bQ f_a(X)$ of $f_a(X)$ over $\bQ$ are totally real cyclic 
quartic number fields which are called the simplest quartic fields 
(cf. e.g. \cite{Gra77}, \cite{Laz91}, \cite{LP95}, 
\cite[Section 6.2]{Gaa02}, \cite{Kim04}, \cite{HH05}, \cite{Duq07}, \cite{Lou07}). 

For $b=-a\in K$, the polynomials $f_a(X)$ and $f_b(X)$ have the same splitting field over $K$. 
In the present paper, we consider the field isomorphism problem of $f_s(X)$, i.e. 
for a fixed $a\in K$, determine whether $b\in K$ gives the same splitting field over $K$ as 
$\mathrm{Spl}_K f_a(X)=\mathrm{Spl}_K f_b(X)$ or not. 

For $n\geq 3$, Rikuna \cite{Rik02} constructed one-parameter families of cyclic 
polynomials of degree $n$ over $K$ with char $K {\not|}$ $n$ and 
$K\ni\zeta+\zeta^{-1}$ where $\zeta$ is a primitive $n$-th root of unity, 
and $f_s(X)$ may be obtained the quartic case $n=4$ of Rikuna's cyclic polynomials 
(see also \cite{Miy99}, \cite{HM99}). 
An answer to the field isomorphism problem to Rikuna's cyclic polynomials 
was given by Komatsu \cite{Kom04} as a generalization of Kummer theory 
(cf. also \cite{Oga03}, \cite{Kid05}). 

In Section \ref{seSimPoly}, by using multi-resolvent polynomials, 
we will give an explicit form of an answer to the field isomorphism problem of $f_s(X)$ 
over $K$ as the special case of the field intersection problem 
(cf. the simplest cubic case \cite{Mor94}, \cite{Cha96}, \cite{Oga03}, \cite{Kom04}, 
\cite{HM09a}, \cite{H}). 
One of the advantages of using multi-resolvent polynomials is the validity for 
non-abelian groups (see \cite{HM07}, \cite{HM09b}, \cite{HM09c}, \cite{HM}). 

\begin{theorem}\label{thC4}
Let $K$ be a field of char $K\neq 2$ and 
$f_a(X)=X^4-aX^3-6X^2+aX+1\in K[X]$ for $a\in K$. 
For $a,b\in K$ with $a\neq\pm b$ and $(a^2+16)(b^2+16)\neq 0$, 
the following three conditions are equivalent$\,:$\\
{\rm (i)} the splitting fields of $f_a(X)$ and of $f_b(X)$ over $K$ coincide$\,;$\\
{\rm (ii)} the polynomial $f_A(X)$ splits completely into four linear factors over $K$ 
for $A=A_1$ or $A=A_2$ where
\[
A_1=\frac{ab+16}{-a+b}\ \ \mathit{and}\ \ A_2=\frac{ab-16}{a+b}\,;
\]
{\rm (iii)} there exists $z\in K$ such that 
\[
B\,=\,a+\frac{(a^2+16)z(z+1)(z-1)}{f_a(z)}
\]
where $B=b$ or $B=-b$. 

Moreover if $\mathrm{Gal}_K f_a(X)\cong C_4$ $($resp. $\mathrm{Gal}_K f_a(X)\cong C_2$ or $\{1\})$ 
then {\rm (ii)} occurs for only one of $A_1$ and $A_2$ $($resp. for both of $A_1$ and $A_2)$ 
and {\rm (iii)} occurs for only one of $b$ and $-b$ $($resp. for both of $b$ and $-b)$. 
\end{theorem}

Note that the equivalence of the conditions (i) and (iii) is valid also for $a=\pm b$. 

By Theorem \ref{thC4}, for a fixed $a\in K$ with $a^2+16\neq 0$, 
we have $\mathrm{Spl}_K f_b(X)=\mathrm{Spl}_K f_a(X)$ where $b$ 
is given as in Theorem \ref{thC4} (iii) for arbitrary $z\in K$ 
with $f_a(z)\neq 0$ and $b^2+16\neq 0$. 
Hence we have the following: 
\begin{corollary}\label{cor1}
Let $K$ be an infinite field of char $K\neq 2$. 
For a fixed $a\in K$ with $a^2+16\neq 0$, there exist infinitely many $b\in K$ 
such that $\mathrm{Spl}_{K} f_b(X)=\mathrm{Spl}_K f_a(X)$. 
\end{corollary}

The following theorem is well-known as Siegel's theorem for curves of genus zero 
(cf. \cite[Theorem 6.1]{Lan78}, \cite[Chapter 8, Section 5]{Lan83}). 
\begin{theorem-nn}[Siegel]
Let $K$ be a number field and $\mathcal{O}_K$ the ring of integers in $K$. 
If a rational function $\varphi(s)\in K(s)$ has at least three distinct poles, 
then there are only finitely many $z\in K$ such that $\varphi(z)\in\mathcal{O}_K$. 
\end{theorem-nn}

In contrast with Corollary \ref{cor1}, by applying Siegel's theorem to Theorem \ref{thC4}, 
we get:
\begin{corollary}\label{cor3}
Let $K$ be a number field and $\mathcal{O}_K$ the ring of integers in $K$. 
Assume that $a\in \mathcal{O}_K$ with $a^2+16\neq 0$. Then 
there exist only finitely many integers $b\in\mathcal{O}_K$ such that 
$\mathrm{Spl}_{K} f_b(X)=\mathrm{Spl}_K f_a(X)$. 
In particular, there exist only finitely many integers $b\in\mathcal{O}_K$ 
such that $f_{A_i}(X)$, $(i=1,2)$, has a linear factor over $\bQ$ where $A_i$ is given in 
Theorem \ref{thC4} {\rm (ii)}. 
\end{corollary}

We treat the case of $K=\bQ$ and $a=m\in\bZ$. 
We get an application of Theorem \ref{thC4} to a related family of 
quartic Thue equations as follows.

Consider the parametric family of quartic Thue equations
\begin{align*}
F_m(X,Y):=X^4-mX^3Y-6X^2Y^2+mXY^3+Y^4=c
\end{align*}
for $m$, $c\in\bZ$ with $c\neq 0$. 
Note that $f_m(X)=F_m(X,1)$. 
The equation $F_m(X,Y)=c$ has the following solutions 
\begin{align*}
F_m(0,\pm e)&=F_m(\pm e,0)=e^4,\qquad 
F_m(\mp e,\pm e)=F_m(\pm e,\pm e)=-4e^4.
\end{align*}
We call such solutions $(x,y)$ to $F_m(x,y)=c$ with $xy(x+y)(x-y)=0$ 
the \textit{trivial} solutions. 

For $c\in\{\pm 1,\pm 4\}$, Lettl-Peth\"o \cite{LP95} and Chen-Voutier \cite{CV97} 
gave a complete solution to Thue equation $F_m(X,Y)=c$ independently (cf. Section \ref{seThue}). 

In \cite{LPV99}, Lettl-Peth\"o-Voutier showed that for $m\geq 58$, the only primitive solutions 
$(x,y)\in\bZ^2$, i.e. $\mathrm{gcd}(x,y)=1$, to the Thue inequality
\begin{align*}
|F_m(x,y)|\leq 6m+7
\end{align*}
with $|x|\leq y$ are trivial solutions $(0,1)$, $(\pm 1,1)$ and $(\pm 1,2)$. 
Note that $F_m(\pm 1,2)=\pm 6m-7$. 
Wakabayashi \cite{Wak07} investigated Thue inequalities 
$|F_{l,m}(x,y)|\leq k$ with two parameters $l,m$ and $F_{1,m}=F_m$. 

For $m\in\bZ$, we put 
\[
L_m:=\mathrm{Spl}_\bQ f_m(X). 
\]  
We give the following correspondence between integer solutions to $F_m(X,Y)=c$ 
and isomorphism classes of the simplest quartic fields $L_m$. 

\begin{theorem}\label{th2}
Let $m\in\bZ\setminus\{0,\pm 3\}$ and $L_m=\mathrm{Spl}_\bQ f_m(X)$. 
There exists an integer $n\in\bZ\setminus\{\pm m\}$ such that 
$L_n=L_m$ if and only if there exists non-trivial solution $(x,y)\in\bZ^2$, 
i.e. $xy(x+y)(x-y)\neq 0$, to the quartic Thue equation 
\begin{align}
F_m(x,y)=c\tag{$*$}
\end{align}
where $c$ is a divisor of $4(m^2+16)$. 
Moreover integers $m$, $n$ and solutions $(x,y)\in \bZ^2$ to $(*)$ can be chosen 
to satisfy the equation
\begin{align}
N=m+\frac{(m^2+16)xy(x+y)(x-y)}{F_m(x,y)}\tag{$**$}
\end{align}
where either $N=n$ or $N=-n$, and the equation $(**)$ occurs for 
only one of $N=n$ and $N=-n$. 
\end{theorem}
The assumption $m\neq 0,\pm 3$ ensures that $f_m(X)$ is irreducible over $\bQ$, 
$\mathrm{Gal}_\bQ f_m(X)\cong C_4$ and the equality $(**)$ holds for only one of 
$N=n$ and $N=-n$. 
This phenomenon comes from the group theoretical reason (see Section \ref{sePre}). 
Indeed, in the case of $m=\pm 3$, 
$f_{\pm 3}(X)=(X^2\pm X-1)(X^2\mp 4X-1)$ and 
the equation $(**)$ occurs for both of $N=3$ and $N=-3$. 
Thus non-trivial solutions $(x,y)\in\bZ^2$ to $(*)$ which satisfy $(**)$ for $N=-m$ exist 
(cf. Theorem \ref{thC4}). 
Hence the assumption $m\neq 0,\pm 3$ also ensures that a non-trivial solution $(x,y)$ 
to $(*)$ corresponds $N\in\bZ\setminus\{\pm m\}$ via $(**)$. 

If there exists an integer $n\in\bZ\setminus\{\pm m\}$ such that $L_n=L_m$, then 
we may choose a primitive solution $(x,y)\in\bZ^2$ to $(*)$ with $(x,y)\equiv (1,0)\pmod{2}$. 
Then four solutions $\pm(x,y)$, $\pm(y,-x)$ to $(*)$ for $c=d$ and four solutions 
$\pm(x',y')$, $\pm(y',-x')$ to $(*)$ for $c=-4d$ are primitive, 
where $(x',y')=(x-y,x+y)$ and $d$ is an odd divisor of $m^2+16$, and only 
these eight primitive solutions satisfy $(**)$ for the same $N$ as in Theorem \ref{th2}. 

\begin{corollary}\label{corN}
For $m\in\bZ\setminus\{0,\pm 3\}$, let $\mathcal{N}$ be the number of primitive solutions 
$(x,y)\in\bZ^2$ with $xy(x+y)(x-y)\neq 0$ to $F_m(x,y)=c$ where $c$ is a divisor of $4(m^2+16)$. 
Then we have 
\begin{align*}
\#\Bigl\{n\in\bZ\setminus\{\pm m\}\ \Big|\ L_n=L_m,\ n>0\Bigr\}=\frac{\mathcal{N}}{8}
\end{align*}
where $L_m=\mathrm{Spl}_\bQ f_m(X)$. 
In particular, if there does not exist $n\in\bZ\setminus\{\pm m\}$ with $L_n=L_m$ then 
$F_m(x,y)=c$ where $c$ is a divisor of $4(m^2+16)$ has only trivial solutions $(x,y)\in\bZ^2$ 
with $xy(x+y)(x-y)=0$. 
\end{corollary}

However we do not know non-trivial solutions to $(*)$ with $m\geq 0$ except for 
$m=1$, $2$, $3$, $4$, $22$, $103$, $956$. 
By Theorem \ref{thC4}, we can check that $L_1=L_{103}$, $L_2=L_{22}$, $L_4=L_{956}$. 
Some numerical examples will be given in Sections \ref{seprim}, \ref{seRed} and \ref{seComp}.
\section{Preliminaries}\label{sePre}

In order to prove Theorem \ref{thC4}, 
we recall known results of the resolvent polynomials which are 
fundamental tools in the computational aspects of Galois theory 
(cf. \cite{Coh93}, \cite{Coh00}, \cite{Ade01}). 
We intend to explain how to get an answer to the field intersection problem of 
$f_s(X)=X^4-sX^3-6X^2+sX+1$, i.e. for $a,b\in K$ how to determine the intersection of 
$\mathrm{Spl}_K f_a(X)$ and $\mathrm{Spl}_K f_b(X)$. 
An answer to the field isomorphism problem (Theorem \ref{thC4}) may be obtained 
as the special case of the field intersection problem. 

Let $\overline{K}$ be a fixed algebraic closure of a field $K$.
Let $f(X):=\prod_{i=1}^m(X-\alpha_i) \in K[X]$ be a separable polynomial of degree $m$ with 
some fixed order of the roots $\alpha_1,\ldots,\alpha_m\in \overline{K}$. 
By resolvent polynomials with suitable invariants, we may determine the Galois group 
of the polynomial $f(X)$ over $K$ as follows. 

Let $R:=K[x_1,\ldots,x_m]$ be the polynomial ring over $K$ with $m$ variables $x_1,\ldots,x_m$. 
For $\Theta \in R$, we take a surjective homomorphism $\omega_f : R \rightarrow 
k(\alpha_1,\ldots,\alpha_m),\, \Theta(x_1,\ldots,x_m)\mapsto 
\Theta(\alpha_1,\ldots,\alpha_m)$, which is called the specialization map. 
The kernel of $\omega_f$ is the ideal 
$I_f:=\{\Theta\in R \mid \Theta(\alpha_1,\ldots,\alpha_m)=0\}$ in $R$. 

Let $S_m$ be the symmetric group of degree $m$. 
We extend the action of $S_m$ on $m$ letters $\{1,\ldots,m\}$ to that on $R$ by 
$\pi(\Theta(x_1,\ldots,x_m)):=\Theta(x_{\pi(1)},\ldots,x_{\pi(m)})$. 
We define the Galois group of $f(X)$ over $K$ by 
$\mathrm{Gal}(f/K):={\{\pi\in S_m \mid \pi(I_f)\subseteq I_f\}}$. 
Then the Galois group of the splitting field $\mathrm{Spl}_K f(X)$ of $f(X)$ over $K$ 
is isomorphic to $\mathrm{Gal}(f/K)$. 
If we take another ordering of roots $\alpha_{\pi(1)},\ldots,\alpha_{\pi(m)}$ 
of $f(X)$ for some $\pi\in S_m$, the corresponding realization of $\mathrm{Gal}(f/K)$ 
is conjugate in $S_m$. 
Hence, for arbitrary ordering of the roots of $f(X)$, $\mathrm{Gal}(f/K)$ 
is determined up to conjugacy in $S_m$. 

For $H\leq U\leq S_m$, an element $\Theta\in R$ is called a $U$-primitive $H$-invariant if 
$H=\mathrm{Stab}_U(\Theta)$ $:=$ $\{\pi\in U\ |\ \pi(\Theta)=\Theta\}$. 
For a $U$-primitive $H$-invariant $\Theta$, the polynomial 
\[
\mathcal{RP}_{\Theta,U}(X)=\prod_{\opi\in U/H}(X-\pi(\Theta))\in R^U[X]
\]
where $\opi$ runs through the left cosets of $H$ in $U$, 
is called the {\it formal} $U$-relative $H$-invariant resolvent by $\Theta$. 
The polynomial 
\[
\mathcal{RP}_{\Theta,U,f}(X):=\omega_f(\mathcal{RP}_{\Theta,U}(X))
\]
is called the $U$-relative $H$-invariant resolvent of $f$ by $\Theta$. 
The following theorem is fundamental in the theory of resolvent polynomials 
(see e.g. \cite[p.95]{Ade01}). 

\begin{theorem}\label{thfun}
Let $G=\mathrm{Gal}(f/K)$, $H\leq U\leq S_m$ be finite groups with 
$G\leq U$ and $\Theta$ a $U$-primitive  $H$-invariant. 
Suppose that $\mathcal{RP}_{\Theta,U,f}(X)=\prod_{i=1}^l h_i^{e_i}(X)$ gives the decomposition 
of $\mathcal{RP}_{\Theta,U,f}(X)$ into a product of powers of distinct irreducible 
polynomials $h_i(X)$, $(i=1,\ldots,l)$, in $K[X]$. 
Then we have a bijection 
\begin{align*}
G\backslash U/H\quad &\longrightarrow \quad \{h_1^{e_1}(X),\ldots,h_l^{e_l}(X)\}\\
G\, \pi\, H\quad &\longmapsto\quad h_\pi(X)
=\prod_{\tau H\subseteq G\,\pi\,H}\bigl(X-\omega_{f}(\tau(\Theta))\bigr)
\end{align*}
where the product runs through the left cosets $\tau H$ of $H$ in $U$ contained in 
$G\, \pi\, H$, that is, through $\tau=\pi_\sigma \pi$ where $\pi_\sigma$ runs 
a system of representative of the left cosets of $G \cap \pi H\pi^{-1};$ each 
$h_\pi(X)$ is irreducible or a power of an irreducible polynomial with $\mathrm{deg}(h_\pi(X))$ 
$=$ $|G\, \pi\, H|/|H|$ $=$ $|G|/|G\cap \pi H\pi^{-1}|$. 
\end{theorem}
\begin{corollary} 
If $G\leq \pi H\pi^{-1}$ for some $\pi\in U$ then $\mathcal{RP}_{\Theta,U,f}(X)$ has 
a linear factor over $K$. 
Conversely, if $\mathcal{RP}_{\Theta,U,f}(X)$ has a non-repeated linear factor over $K$ 
then there exists $\pi\in U$ such that $G\leq \pi H\pi^{-1}$. 
\end{corollary}

Note that when $\mathcal{RP}_{\Theta,U,f}(X)$ is not squarefree, there exists 
a suitable Tschirnhausen transformation $\hat{f}$ of $f$ over $K$ such that 
$\mathcal{RP}_{\Theta,U,\hat{f}}(X)$ is squarefree (cf. \cite{Gir83}, 
\cite[Alg. 6.3.4]{Coh93}). \\

We apply Theorem \ref{thfun} to the cyclic quartic case. 
Let $f^1(X)\in K[X]$ and $f^2(X)\in K[X]$ be separable quartic polynomials over $K$ respectively. 

We put 
\[
f(X):=f^1(X)f^2(X)
\]
and 
\begin{align*}
G_1:=\mathrm{Gal}(f^1/K),\quad G_2:=\mathrm{Gal}(f^2/K),\quad G:=\mathrm{Gal}(f/K).
\end{align*}

We assume that $G_1,G_2\leq C_4$ 
and apply Theorem \ref{thfun} to $m=8$, $f(X)=f^1(X)f^2(X)$, 
$U=\langle\sigma\rangle\times\langle\tau\rangle$, $H=\langle\sigma\tau\rangle$ or 
$\langle\sigma\tau^3\rangle$ where $\sigma,\tau\in S_8$ act on $R=K[x_1,\ldots,x_8]$ by 
\begin{align*}
\sigma\,:\, x_1\mapsto x_2\mapsto x_3\mapsto x_4\mapsto x_1,\\
\tau\,:\, x_5\mapsto x_6\mapsto x_7\mapsto x_8\mapsto x_5.
\end{align*}
We put $U:=\langle\sigma\rangle\times\langle\tau\rangle$. 
Let $\Theta_1$ (resp. $\Theta_2$) be a $U$-primitive 
$\langle\sigma\tau\rangle$-invariant (resp. $\langle\sigma\tau^3\rangle$-invariant). 
Then we have the $U$-relative $\langle\sigma\tau\rangle$-invariant 
(resp. $\langle\sigma\tau^3\rangle$-invariant) resolvent 
polynomial of $f(X)=f^1(X)f^2(X)$ by $\Theta_1$ (resp. $\Theta_2$) as 
\begin{align*}
\mathcal{R}_f^i(X):=
\mathcal{RP}_{\Theta_i,U,f}(X),\quad 
(i=1,2).
\end{align*}
This kind of resolvent polynomial is also called (absolute) 
{\it multi-resolvent polynomial} (cf. \cite{RV99}, \cite{Ren04}). 

For a squarefree polynomial $\mathcal{R}(X)\in K[X]$ of degree $l$, we define the 
{\it decomposition type} $\mathrm{DT}(\mathcal{R})$ of $\mathcal{R}(X)$ by the partition 
of $l$ induced by the degrees of the irreducible factors of $\mathcal{R}(X)$ over $K$. 
By Theorem \ref{thfun}, we get the intersection field $\mathrm{Spl}_K f^1(X)$ $\cap$ 
$\mathrm{Spl}_Kf^2(X)$ via the decomposition types $\mathrm{DT}(\mathcal{R}_f^1)$ 
and $\mathrm{DT}(\mathcal{R}_f^2)$. 
\begin{theorem}\label{thDecom}
For $f(X)=f^1(X)f^2(X)\in K[X]$ with $G_1$, $G_2\leq C_4$, 
we assume that $\#G_1\geq\#G_2$ and both $\mathcal{R}_f^1(X)$ and 
$\mathcal{R}_f^2(X)$ are squarefree. 
Then the Galois group $G=\mathrm{Gal}(f/K)$ and the intersection field 
$\mathrm{Spl}_K f^1(X)\cap\mathrm{Spl}_K f^2(X)$ are given by the decomposition 
types $\mathrm{DT}(\mathcal{R}_f^1)$ and $\mathrm{DT}(\mathcal{R}_f^2)$ 
as on Table $1$. 
\end{theorem}
\begin{center}
{\rm Table} $1$\vspace*{5mm}\\
{\renewcommand\arraystretch{1.1}
\begin{tabular}{|c|c|c|l|l|l|}\hline
$G_1$& $G_2$ & $G$ & & ${\rm DT}(\mathcal{R}_f^1)$ 
& ${\rm DT}(\mathcal{R}_f^2)$ \\ \hline 
& & $C_4\times C_4$ & $L_1\cap L_2=K$ & $4$ & $4$\\\cline{3-6} 
& \raisebox{-1.6ex}[0cm][0cm]{$C_4$} & $C_4\times C_2$ & $[L_1\cap L_2:K]=2$ 
& $2,2$ & $2,2$\\ \cline{3-6} 
& & \raisebox{-1.6ex}[0cm][0cm]{$C_4$} 
& \raisebox{-1.6ex}[0cm][0cm]{$L_1=L_2$} & $2,2$ & $1,1,1,1$\\ \cline{5-6}
$C_4$ & & & & $1,1,1,1$ & $2,2$\\ \cline{2-6} 
& \raisebox{-1.6ex}[0cm][0cm]{$C_2$}  & $C_4\times C_2$ & $L_1\cap L_2=K$ 
& $4$ & $4$\\ \cline{3-6}
& & $C_4$ & $L_1\supset L_2$ & $4$ & $4$\\ \cline{2-6} 
& $\{1\}$  & $C_4$ & $L_1\supset L_2=K$ & $4$ & $4$\\ \cline{1-6}
& \raisebox{-1.6ex}[0cm][0cm]{$C_2$}  & $C_2\times C_2$ & $L_1\cap L_2=K$ 
& $2,2$ & $2,2$\\ \cline{3-6}
$C_2$ & & $C_2$ & $L_1=L_2$ & $1,1,1,1$ & $1,1,1,1$\\ \cline{2-6}
& \{1\} & $C_2$ & $L_1\supset L_2$ & $2,2$ & $2,2$\\ \cline{1-6}
$\{1\}$ & $\{1\}$ & $\{1\}$ & $L_1=L_2=K$ & $1,1,1,1$ & $1,1,1,1$\\ \cline{1-6}
\end{tabular}
}
\vspace*{4mm}
\end{center}

We checked the decomposition types $\mathrm{DT}(\mathcal{R}_f^i)$, $(i=1,2)$, 
on Table $1$ by GAP \cite{GAP}. 

Now we get an answer to the field isomorphism problem of 
\[
f_s(X)=X^4-sX^3-6X^2+sX+1
\]
via multi-resolvent polynomials $\mathcal{R}_{f_{a,b}}^i(X):=
\mathcal{RP}_{\Theta_i,\langle\sigma\rangle\times\langle\tau\rangle,f_{a,b}}$, 
($i=1,2$), where 
\[
f_{a,b}(X):=f_a(X)f_b(X)
\] 
as the special case of Theorem \ref{thDecom}. 
Note that $\mathrm{disc}(f_s(X))=4(s^2+16)^3$. 

\begin{theorem}\label{thRab}
For $a,b\in K$ with $(a^2+16)(b^2+16)\neq 0$, 
we assume that both $\mathcal{R}_{f_{a,b}}^1(X)$ and $\mathcal{R}_{f_{a,b}}^2(X)$ 
are squarefree. 
Then two splitting fields of $f_a(X)$ and of $f_b(X)$ over $K$ coincide if and only if 
$\mathcal{R}_{f_{a,b}}^1(X)$ or $\mathcal{R}_{f_{a,b}}^2(X)$ splits completely into four 
linear factors over $K$. 
\end{theorem}

This is an analogue of the classical result by Kummer theory. 
Namely for a field $K$ which contains a primitive fourth root of unity and $a,b\in K$, 
the splitting fields of $X^4-a$ and of $X^4-b$ over $K$ coincide 
if and only if $X^4-ab$ or $X^4-ab^3$ has a linear factor 
(equivalent to split completely) over $K$. 
It is remarkable that Theorem \ref{thRab} does not need the assumption that 
$K$ contains a primitive fourth root of unity. 

\section{Proof of Theorem \ref{thC4}}\label{seSimPoly}

We give an explicit answer to the field intersection problem of the simplest quartic 
polynomials $f_s(X)$ via suitable invariants $\Theta_1$ and $\Theta_2$. 
As the special case, we obtain Theorem \ref{thC4}. 

Let $K(z)$ be the rational function field over $K$ 
and $\sigma$ a $K$-automorphism of $K(z)$ of order four which is defined by 
\[
\sigma : z\ \mapsto\ \frac{z-1}{z+1}\ \mapsto\ -\frac{1}{z}\ 
\mapsto\ -\frac{z+1}{z-1}\ \mapsto z.
\]
We consider the fixed field $K(z)^{\langle\sigma\rangle}$ 
and the $C_4$-extension $K(z)/K(z)^{\langle\sigma\rangle}$. 
Then we get 
\begin{align*}
f_s(X)&=\prod_{x\in\mathrm{Orb}_{\langle\sigma\rangle}(z)}\Bigl(X-x\Bigr)
=\Bigl(X-z\Bigr)\Bigl(X-\frac{z-1}{z+1}\Bigr)\Bigl(X+\frac{1}{z}\Bigr)
\Bigl(X+\frac{z+1}{z-1}\Bigr)\\
&=X^4-sX^3-6X^2+sX+1
\end{align*}
where
\begin{align*}
s=\frac{z^4-6z^2+1}{z(z^2-1)}=\frac{(z^2+2z-1)(z^2-2z-1)}{z(z+1)(z-1)}
\end{align*}
as the generating polynomial of the field extension $K(z)/K(z)^{\langle\sigma\rangle}$. 
It follows that $K(z)^{\langle\sigma\rangle}=K(s)$ and the Galois group of the polynomial 
$f_s(X)$ over $K(s)$ is isomorphic to $C_4$. 

We also take another rational function field $K(w)$ over $K$ with indeterminate $w$, 
$\tau\in\mathrm{Aut}_K K(w)$ with 
\begin{align*}
\tau : w\ \mapsto\ \frac{w-1}{w+1}\ \mapsto\ -\frac{1}{w}\ 
\mapsto\ -\frac{w+1}{w-1}\ \mapsto w
\end{align*}
and $f_t(X)=X^4-tX^3-6X^2+tX+1$ where 
\begin{align*}
t=\frac{w^4-6w^2+1}{w(w^2-1)}=\frac{(w^2+2w-1)(w^2-2w-1)}{w(w+1)(w-1)}
\end{align*}
by the same manner of $K(z)$, $\sigma$ and $f_s(X)$. 

Put $U:=\langle\sigma\rangle\times\langle\tau\rangle$. 
Then the field $K(z,w)$ is $(C_4\times C_4)$-extension of $K(z,w)^U=K(s,t)$. 

In order to apply Theorem \ref{thRab}, we should find suitable $U$-primitive 
$\langle\sigma\tau\rangle$-invariant $\Theta_1$ and 
$U$-primitive $\langle\sigma\tau^3\rangle$-invariant $\Theta_2$.  

We may find the following two $U$-primitive $\langle\sigma\tau\rangle$-invariants 
which are candidates to $\Theta_1$: 
\begin{align*}
\Theta_1&=\sum_{i=0}^3(\sigma\tau)^i(zw)
=\frac{(w+z)(wz-1)(zw-w-z-1)(zw+w+z-1)}{zw(z^2-1)(w^2-1)},\quad\mathrm{or}\\
\Theta_1&=\prod_{i=0}^3(\sigma\tau)^i(z+w)\\
&=\frac{(z^2w^2-zw^2-z^2w+2zw+w+z+1)(z^2w^2+zw^2+z^2w+2zw-w-z+1)}{zw(z^2-1)(w^2-1)}.
\end{align*}

However the multi-resolvent polynomial $\mathcal{R}_{f_{a,b}}^1(X)
=\mathcal{RP}_{\Theta_1,U,f_{a,b}}(X)$ where $f_{a,b}(X)=f_a(X)f_b(X)$ 
becomes complicated in the both cases above. 

\begin{example}\label{ex1}

We present two explicit examples of the multi-resolvent polynomials 
$\mathcal{R}_{f_{a,b}}^i(X)$ $:=$ 
$\mathcal{RP}_{\Theta_i,\langle\sigma\rangle\times\langle\tau\rangle,f_{a,b}}(X)$ for 
$i=1,2$ where $f_{a,b}(X):=f_a(X)f_b(X)$.\\

(i) In \cite{HM09b}, we gave an answer to the field isomorphism problem of $f_s(X)$ 
by taking 
\begin{align*}
\Theta_1=\frac{2(z^2+1)(w^2+1)(zw+1)(z-w)}{zw(z^2-1)(w^2-1)},\quad 
\Theta_2=\frac{2(z^2+1)(w^2+1)(zw-1)(z+w)}{zw(z^2-1)(w^2-1)}.
\end{align*}
Then the corresponding multi-resolvent polynomials are given as 
\begin{align*}
\mathcal{R}_{f_{a,b}}^1(X)=X^4-(a^2+16)(b^2+16)(X^2-4(a-b)^2),\\
\mathcal{R}_{f_{a,b}}^2(X)=X^4-(a^2+16)(b^2+16)(X^2-4(a+b)^2).
\end{align*}

(ii) If we take another $U$-primitive $\langle\sigma\tau\rangle$-invariant $\Theta_1$ 
and $\langle\sigma\tau^3\rangle$-invariant $\Theta_2$ as 
\begin{align*}
\Theta_1=\frac{2(zw+1)(z-w)}{(z^2+1)(w^2+1)},\quad 
\Theta_2=\frac{2(zw-1)(z+w)}{(z^2+1)(w^2+1)}
\end{align*}
then we get 
\begin{align*}
\mathcal{R}_{f_{a,b}}^1(X)=X^4-X^2+\frac{4(a-b)^2}{(a^2+16)(b^2+16)},\quad 
\mathcal{R}_{f_{a,b}}^2(X)=X^4-X^2+\frac{4(a+b)^2}{(a^2+16)(b^2+16)}.
\end{align*}
\end{example}

The multi-resolvent polynomials in Example \ref{ex1} are useful since they are biquadratic, 
i.e. quadratic polynomial with respect to $X^2$. 
However we do not understand for a fixed $a\in K$ whether there exist infinitely many 
$b\in K$ such that $\mathcal{R}_{f_{a,b}}^i(X)$, ($i=1,2$), splits completely over $K$ or not. 
By Theorem \ref{thRab}, this question means that 
for a fixed $a\in K$ whether there exist infinitely many 
$b\in K$ such that $\mathrm{Spl}_K f_a(X)=\mathrm{Spl}_K f_b(X)$.\\

It follows from \cite[Theorem 1.4]{AHK98} that there exist 
$\langle\sigma\tau\rangle$-invariant $\Theta_1$ and 
$\langle\sigma\tau^3\rangle$-invariant $\Theta_2$ such that 
$K(z,w)=K(z,\Theta_1)=K(z,\Theta_2)$. 
The following gives such invariants $\Theta_1$ and $\Theta_2$ 
which is the key lemma of this paper. 
\begin{lemma}\label{lemTheta}
Let $U=\langle\sigma\rangle\times\langle\tau\rangle$. 
We take 
\[
\Theta_1:=\frac{zw+1}{-z+w}\quad\text{and}\quad 
\Theta_2:=\frac{zw-1}{z+w}.
\]
Then the following assertions hold$\,:$\\
{\rm (i)} the element $\Theta_1$ is a $U$-primitive $\langle\sigma\tau\rangle$-invariant$\,;$\\
{\rm (ii)} the element $\Theta_2$ is a $U$-primitive $\langle\sigma\tau^3\rangle$-invariant$\,;$\\
{\rm (iii)} the $U$-orbit of $\Theta_i$ is given by 
the same as $\langle\sigma\rangle$-orbit of $z\,;$ 
\begin{align*}
\mathrm{Orb}_U(\Theta_i)=
\Big\{\Theta_i,\ \frac{\Theta_i-1}{\Theta_i+1},\ -\frac{1}{\Theta_i},\ 
-\frac{\Theta_i+1}{\Theta_i-1}\Big\},\quad (i=1,2). 
\end{align*}
\end{lemma}
\begin{proof}
We can check the assertions by direct computations. 
\end{proof}

The multi-resolvent polynomials $\mathcal{R}_{f_{a,b}}^i(X)$ 
$:=$ $\mathcal{RP}_{\Theta_i,\langle\sigma\rangle\times\langle\tau\rangle,f_af_b}(X)$, 
$(i=1,2)$, with respect to
\[
\Theta_1=\frac{zw+1}{-z+w}\quad\text{and}\quad 
\Theta_2=\frac{zw-1}{z+w}
\] 
as in Lemma \ref{lemTheta} are given by 
\begin{align}
\mathcal{R}_{f_{a,b}}^1(X)&=f_{A_1}(X)=X^4-\frac{ab+16}{-a+b}X^3
-6X^2+\frac{ab+16}{-a+b}X+1,\label{polyR}\\
\mathcal{R}_{f_{a,b}}^2(X)&=f_{A_2}(X)=X^4-\frac{ab-16}{a+b}X^3
-6X^2+\frac{ab-16}{a+b}X+1\nonumber
\end{align}
where 
\[
A_1=\frac{ab+16}{-a+b}\ \ \mathit{and}\ \ A_2=\frac{ab-16}{a+b}. 
\]
Note that 
\[
\mathrm{disc}(\mathcal{R}_{f_{a,b}}^1(X))=\frac{4(a^2+16)^3(b^2+16)^3}{(a-b)^6},\quad 
\mathrm{disc}(\mathcal{R}_{f_{a,b}}^2(X))=\frac{4(a^2+16)^3(b^2+16)^3}{(a+b)^6}.
\]

By Theorem \ref{thDecom}, we get the intersection field $\mathrm{Spl}_K f_a(X)$ 
$\cap$ $\mathrm{Spl}_K f_b(X)$ via Table $1$. 

\begin{theorem}\label{thDab}
Let $\mathcal{R}_{f_{a,b}}^i(X)$, $(i=1,2)$, be as in $(\ref{polyR})$.
For $a,b\in K$ with $a\neq\pm b$ and $(a^2+16)(b^2+16)\neq 0$, 
we assume that $\#\mathrm{Gal}_K f_a(X)\geq\#\mathrm{Gal}_K f_b(X)$. 
Then the intersection field 
$\mathrm{Spl}_K f_a(X)\cap\mathrm{Spl}_K f_b(X)$ is given by the decomposition 
types $\mathrm{DT}(\mathcal{R}_{f_{a,b}}^1)$ and $\mathrm{DT}(\mathcal{R}_{f_{a,b}}^2)$ 
as on Table $1$ in Theorem \ref{thDecom}.
\end{theorem}

\begin{proof}[Proof of Theorem \ref{thC4}]
As the special case of Theorem \ref{thDab}, we see the conditions (i) and (ii) are equivalent 
(cf. also Theorem \ref{thRab}). 

The condition (iii) is just a restatement of (ii). 
Indeed, we may check that 
$z\in K$ is a root of $f_{A_1}(X)$ (resp. $f_{A_2}(X)$) if and only if 
$z$ satisfies the condition (iii) for $B=b$ (resp. $B=-b$). 
Note that if $z$ is a root of $f_{A_i}(X)$ then 
$\frac{z-1}{z+1}$, $-\frac{1}{z}$, $-\frac{z+1}{z-1}$ 
are also roots of $f_{A_i}(X)$ for $i=1,2$. 
By Table $1$ as in Theorem \ref{thDecom}, 
if $\mathrm{Spl}_K f_a(X)\cong C_4$ (resp. $C_2$ or \{1\}) and $f_{A_i}(z)$ splits completely 
then $f_{A_j}(X)$ is irreducible (resp. splits completely) over $K$ for 
$(i,j)=(1,2)$ and $(2,1)$. 
This completes the proof.
\end{proof}

\section{Simplest quartic fields and related Thue equations}\label{seThue}

In this section, we treat the case of $K=\bQ$ and $a=m\in\bZ$ for 
\begin{align*}
f_a(X)=f_m(X)=X^4-mX^3-6X^2+mX+1\in\bZ[X].
\end{align*}

We first see
\begin{lemma}
{\rm (i)} For $m\in\bZ\setminus\{0,\pm 3\}$, $f_m(X)$ is irreducible over $\bQ\,;$\\
{\rm (ii)} $f_{\pm 3}(X)=(X^2\pm X-1)(X^2\mp 4X-1)$, $f_0(X)=(X^2+2X-1)(X^2-2X-1)$. 
\end{lemma}
\begin{proof}
It follows from $\mathrm{Gal}_{\bQ(s)} f_s(X)\cong C_4$ that 
$\mathrm{Gal}_\bQ f_m(X)\cong C_4$, $C_2$ or $\{1\}$ for $m\in\bZ$. 
Thus if $f_m(X)$ is reducible over $\bQ$ then there exist $a,b,c\in\bZ$ such that 
$f_m(X)=(X^2+aX+c)(X^2+bX+c)$ with $c=\pm 1$. 
By comparing the coefficients, we see $c=-1$, because if $c=1$ then 
$(m,a,b)=(0,\pm 2\sqrt{2},\mp 2\sqrt{2})$. 
Also by comparing the coefficients, if $c=-1$ then $(m,a)=(-b+4/b,-4/b)\in\bZ^2$. 
Hence $b\in\{\pm 1,\pm 2,\pm 4\}$. 
In each case, we have $(m,a,b)=(3,-4,1)$, $(-3,4,-1)$, $(0,-2,2)$, $(0,2,-2)$, 
$(-3,-1,4)$, $(3,1,-4)$.
\end{proof}
For $m\in\bZ\setminus\{0,\pm 3\}$, the splitting fields 
$L_m:=\mathrm{Spl}_\bQ f_m(X)$ of $f_m(X)$ over $\bQ$ are totally real cyclic 
quartic number fields and called the simplest quartic fields. 

We consider the related parametric family of quartic Thue equations
\begin{align*}
F_m(X,Y):=X^4-mX^3Y-6X^2Y^2+mXY^3+Y^4=c
\end{align*}
for $m$, $c\in\bZ$ with $c\neq 0$. 
Note that $f_m(X)=F_m(X,1)$. 

We may assume that $m\geq 0$ because if $(x,y)\in\bZ^2$ is a solution 
to $F_m(x,y)=c$ then $(y,x)$ becomes a solution to $F_{-m}(y,x)=c$. 
The equation $F_m(x,y)=c$ has the following solutions 
\begin{align*}
F_m(0,\pm e)&=F_m(\pm e,0)=e^4,\qquad 
F_m(\mp e,\pm e)=F_m(\pm e,\pm e)=-4e^4.
\end{align*}
We call such solutions $(x,y)\in\bZ^2$ to $F_m(x,y)=c$ with $xy(x+y)(x-y)=0$ 
the {\it trivial} solutions. 

If $(x,y)\in\bZ^2$ is a solution to $F_m(x,y)=c$ then four pairs 
$\pm(x,y)$, $\pm(y,-x)$ are also solutions because $F_m(X,Y)$ 
is invariant under the action $X\longmapsto Y\longmapsto -X$ of order four. 

For $c\in\{\pm 1,\pm 4\}$, Lettl-Peth\"o \cite{LP95} and Chen-Voutier \cite{CV97} 
gave a complete solution to Thue equation $F_m(X,Y)=c$ independently. 
For $c\in\{\pm 1,\pm 4\}$ and $m\geq 0$, all solutions to $F_m(X,Y)=c$ are given by 
eight trivial solutions 
\begin{align*}
F_m(0,\pm 1)=F_m(\pm 1,0)=1,\qquad F_m(\mp 1,\pm 1)=F_m(\pm 1,\pm 1)=-4
\end{align*}
for arbitrary $m\geq 0$ and additionally 
\begin{align*}
F_1(\mp 2,\pm 1)&=F_1(\pm 1,\pm 2)=-1,& 
F_1(\pm 3,\pm 1)&=F_1(\mp 1,\pm 3)=4,\\
F_4(\mp 3,\pm 2)&=F_4(\pm 2,\pm 3)=1,& 
F_4(\pm 5,\pm 1)&=F_4(\mp 1,\pm 5)=-4.
\end{align*}
We put 
\[
(x',y'):=(x-y,x+y).
\]
Then if $(x,y)\in\bZ^2$ is a solution to $F_m(x,y)=c$ then $(x',y')\in\bZ^2$ is a solution 
to $F_m(x',y')=-4c$. 
Conversely if $(x',y')$ is a solution to $F_m(x',y')=-4c$ then $x'\equiv y'\pmod{2}$ and 
$(x,y)=(\frac{x'+y'}{2},\frac{-x'+y'}{2})\in\bZ^2$ is a solution to $F_m(x,y)=c$. 

In \cite{LPV99}, Lettl-Peth\"o-Voutier showed that for $m\geq 58$, 
the only primitive solutions $(x,y)\in\bZ^2$, i.e. $\mathrm{gcd}(x,y)=1$, 
to the Thue inequality 
\begin{align*}
|F_m(x,y)|\leq 6m+7
\end{align*}
are trivial solutions $\pm (0,1)$, $\pm(1,0)$, $\pm (1,1)$, $\pm(-1,1)$ and 
non-trivial solutions $\pm(2,1)$, $\pm(-1,2)$, $\pm(-2,1)$, $\pm(1,2)$. 
We note that $F_m(\pm 2,1)=F_m(\mp 1,2)=\mp 6m-7$. 

If $(x,y)\in\bZ^2$ is a primitive solution to $(*)$ then four pairs 
$\pm(x,y)$, $\pm(y,-x)$ are primitive solutions to $(*)$. 
We also see the following lemma: 
\begin{lemma}\label{lemprim}
Put $(x',y'):=(x-y,x+y)$. 
If $(x,y)\in\bZ^2$ with $(x,y)\equiv (0,1)$ or $(1,0)\pmod{2}$ is a primitive solution 
to $F_m(x,y)=c$ then $c$ is an odd integer and $(x',y')\equiv (1,1)\pmod{2}$ is 
a primitive solution to $F_m(x',y')=-4c$. 
Conversely if $(x',y')\in\bZ^2$ with $(x',y')\equiv (1,1)\pmod{2}$ is a primitive solution 
to $F_m(x',y')=d$ then $d=-4c$ for an odd integer $c$ and 
$(x,y)=(\frac{x'+y'}{2},\frac{-x'+y'}{2})\equiv (0,1)$ or $(1,0)\pmod{2}$ is 
a primitive solution to $F_m(x,y)=c$. 
\end{lemma}
\begin{proof}
Assume that $(x,y)\in\bZ^2$ with $(x,y)\equiv (0,1)$ or $(1,0)\pmod{2}$ 
is a primitive solution to $F_m(x,y)=c$. 
Then $c$ is odd because $F_m(0,1)=F_m(1,0)=1$. 
It follows by the definition that $(x',y')\equiv(1,1)\pmod{2}$ and $F_m(x',y')=-4c$. 
If $\gcd(x',y')=k'>1$ then $k'$ is odd and $k'$ divides both $x=\frac{x'+y'}{2}$ and 
$y=\frac{-x'+y'}{2}$. 
This contradicts to $\gcd(x,y)=1$. 
Hence we have $\gcd(x',y')=1$. 

Conversely we assume that 
$(x',y')\in\bZ^2$ with $(x',y')\equiv (1,1)\pmod{2}$ is a primitive solution 
to $F_m(x',y')=d$. 
Then $F_m(1,1)=-4\equiv 0\pmod{4}$ and $\not\equiv 0\pmod{8}$. 
Thus $d=-4c$ for an odd integer $c\in\bZ$. 
We also see $F_m(x,y)=c$. 
If $\gcd(x,y)=k>1$ then $k$ divides both $x'=x-y$ and $y'=x+y$. 
This contradicts to $\gcd(x',y')=1$. 
Therefore we have $\gcd(x,y)=1$. 
Because if $(x,y)\equiv(1,1)\pmod{2}$ then $F_m(x,y)\equiv 0\pmod{4}$, 
we also obtain $(x,y)\equiv (0,1)$ or $(1,0)\pmod{2}$. 
\end{proof}

\section{Proof of Theorem \ref{th2}: the correspondence}\label{secorr}

The aim of this section is to establish the correspondence between isomorphism classes 
of the simplest quartic fields $L_m$ and non-trivial solutions to quartic 
Thue equations $(*)$ as follows: 

\begin{theorem-nn}[Theorem \ref{th2}]
Let $m\in\bZ\setminus\{0,\pm 3\}$ and $L_m=\mathrm{Spl}_\bQ f_m(X)$. 
Then there exists an integer $n\in\bZ\setminus\{\pm m\}$ such that 
$L_m=L_n$ if and only if there exists non-trivial solution $(x,y)\in\bZ^2$, 
$i.e.\ xy(x+y)(x-y)\neq 0$, to the quartic Thue equation 
\begin{align}
F_m(x,y)=x^4-mx^3y-6x^2y^2+mxy^3+y^4=c\tag{$*$}
\end{align}
where $c$ is a divisor of $4(m^2+16)$. 
Moreover integers $m$, $n$ and solutions $(x,y)\in \bZ^2$ to $(*)$ can be chosen 
to satisfy the equation
\begin{align}
N=m+\frac{(m^2+16)xy(x+y)(x-y)}{F_m(x,y)}\tag{$**$}
\end{align}
where either $N=n$ or $N=-n$, and the equation $(**)$ occurs for 
only one of $N=n$ and $N=-n$. 
\end{theorem-nn}
%

%
%
%
%
\begin{proof}[Proof of Theorem \ref{th2}]

We use Theorem \ref{thC4} in the case where $K=\bQ$.\\

For $m\in\bZ\setminus\{0,\pm 3\}$, we assume that 
there exists an integer $n\in\bZ\setminus\{\pm m\}$ such that $L_m=L_n$. 
By Theorem \ref{thC4}, there exists $z\in\bQ$ such that 
\[
N\,=\,m+\frac{(m^2+16)z(z+1)(z-1)}{f_m(z)}
\]
where either $N=n$ or $N=-n$. 
Write $z=x/y$ with  $x,y\in\bZ$ and $\mathrm{gcd}(x,y)=1$ then we have
\[
N\,=\,m+\frac{(m^2+16)xy(x+y)(x-y)}{F_m(x,y)}\in\bZ. 
\]
By the assumption $n\neq \pm m$, we have $xy(x+y)(x-y)\neq 0$.\\

We will show that $c:=F_m(x,y)$ divides $4(m^2+16)$. 
We make use of a resultant and the Sylvester matrix (cf. \cite{PV00}, 
\cite[Section 1.3]{SWP08}, see also \cite[Theorem 6.1]{Lan78}, 
\cite[Chapter 8, Section 5]{Lan83}). \\

Put $h(z):=(m^2+16)z(z+1)(z-1)$ and $f(z):=F_m(z,1)$. 
We take the resultant
\[
R_m:=\mathrm{Res}_z(h(z),f(z))=16(m^2+16)^4
\]
of $h(z)$ and $f(z)$ with respect to $z$. 
We see that $R_m$ is also given by 
\begin{align*}
R_m&=
\left|
\begin{array}{lllllll}
 m^2+16 & 0 & -m^2-16 & 0 & 0 & 0 & h(z) z^3 \\
 0 & m^2+16 & 0 & -m^2-16 & 0 & 0 & h(z) z^2 \\
 0 & 0 & m^2+16 & 0 & -m^2-16 & 0 & h(z) z \\
 0 & 0 & 0 & m^2+16 & 0 & -m^2-16 & h(z) \\
 1 & -m & -6 & m & 1 & 0 & f(z) z^2 \\
 0 & 1 & -m & -6 & m & 1 & f(z) z \\
 0 & 0 & 1 & -m & -6 & m & f(z)
\end{array}
\right|\\
&=4(m^2+16)^3\Bigl(h(z)p(z)+f(z)q(z)\Bigr)
\end{align*}
where
\begin{align*}
p(z)=5z^3-5mz^2-29z+4m,\quad q(z)=-(m^2+16)(5z^2-4).
\end{align*}
Hence we have 
\begin{align*}
h(z)p(z)+f(z)q(z)=4(m^2+16).
\end{align*}

Put 
\begin{align*}
H(x,y)&:=(m^2+16)xy(x+y)(x-y),\\
P(x,y)&:=5x^3-5mx^2y-29xy^2+4my^3,\\
Q(x,y)&:=-(m^2+16)y(5x^2-4y^2). 
\end{align*}
Then it follows from $z=x/y$ that 
\begin{align*}
H(x,y)P(x,y)+F_m(x,y)Q(x,y)=4(m^2+16)y^7.
\end{align*}
Because the cubic forms $F_m(X,Y)$ and $H(X,Y)$ are invariants under the action of 
$\sigma\,:\, X\mapsto Y$, $Y\mapsto -X$, we also get 
\begin{align*}
H(x,y)P(y,-x)+F_m(x,y)Q(y,-x)=4(m^2+16)(-x)^7.
\end{align*}

Hence we have
\begin{align*}
&\frac{H(x,y)P(x,y)}{F_m(x,y)}+Q(x,y)=\frac{4(m^2+16)y^7}{F_m(x,y)}\in\bZ,\\
&\frac{H(x,y)P(y,-x)}{F_m(x,y)}+Q(y,-x)=-\frac{4(m^2+16)x^7}{F_m(x,y)}\in\bZ.
\end{align*}

Since $x$ and $y$ are relatively prime, 
we conclude that $c=F_m(x,y)$ divides $4(m^2+16)$. \\

Conversely if there exists $(x,y)\in\bZ^2$ with $xy(x+y)(x-y)\neq 0$ such that 
$c=F_m(x,y)$ divides $4(m^2+16)$ then we can take $N\in\bQ\setminus\{m\}$ by 
\begin{align*}
N\,=\,m+\frac{(m^2+16)xy(x+y)(x-y)}{F_m(x,y)}. 
\end{align*}
From the assumption $m\in\bZ\setminus\{0,\pm 3\}$, we have $\mathrm{Gal}_\bQ f_m(X)\cong C_4$. 
Hence it follows from Theorem \ref{thC4} (i) and (iii) that 
$N\in\bQ\setminus\{\pm m\}$ and $L_m=L_N$. 

We should show that $N\in\bZ$. 
If $x\not\equiv y\pmod{2}$ then $c=F_m(x,y)$ divides $m^2+16$ and hence 
$N\in\bZ\setminus\{\pm m\}$, because $F_m(x,y)\equiv F_m(x,y)\equiv 1\pmod{2}$. 

If $x\equiv y\pmod{2}$ then $c=F_m(x,y)$ divides $(m^2+16)xy(x+y)(x-y)$ 
and hence $N\in\bZ\setminus\{\pm m\}$, because $xy(x+y)(x-y)\equiv 0\pmod{4}$. 
\end{proof}

%
%
\section{Primitive solutions}\label{seprim}

By the proof of Theorem \ref{thC4} and Theorem \ref{th2}, 
a non-trivial solution $(x,y)\in\bZ^2$ to $(*)$ may be obtained as 
$z=x/y$ with $\mathrm{gcd}(x,y)=1$ 
where $z\in\bQ$ is a root of $f_A(X)=X^4-AX^3-6X^2+AX+1$ for $A=A_1$ or $A=A_2$ with 
\[
A_1=\frac{mn+16}{-m+n}\ \ \mathrm{and}\ \ A_2=\frac{mn-16}{m+n}
\]
as in Theorem \ref{thC4} (ii). 
Hence we now consider only primitive solutions $(x,y)\in\bZ^2$, i.e. 
$\mathrm{gcd}(x,y)=1$, to $(*)$.

\begin{lemma}
Let $m\in\bZ\setminus\{0,\pm 3\}$ and $L_m=\mathrm{Spl}_\bQ f_m(X)$. 
We assume that there exists $n\in\bZ\setminus\{\pm m\}$ such that $L_n=L_m$.\\
{\rm (i)} We may choose non-trivial primitive solution $(x,y)\in\bZ^2$ with 
$(x,y)\equiv(0,1)\pmod{2}$ to $F_m(x,y)=d$ where $d$ is an odd divisor of $m^2+16$. 
Then four pairs $\pm(x,y)$, $\pm(y,-x)$ are primitive solutions to $F_m(X,Y)=d$ 
and four pairs $\pm(x',y')$, $\pm(y',-x')$ are primitive solutions to $F_m(X,Y)=-4d$ where 
$(x',y')=(x-y,x+y)$.\\
{\rm (ii)} All primitive solutions to $(*)$ which satisfy $(**)$ for either $N=n$ or 
$N=-n$ are given by the eight solutions as in {\rm (i)}, and such solutions exist for 
only one of $N=n$ and $N=-n$. 
\end{lemma}
\begin{proof}
By Theorem \ref{th2}, there exists non-trivial solution $(x_0,y_0)\in\bZ^2$ to $F_m(x_0,y_0)=c$ 
for a divisor $c$ of $4(m^2+16)$. 
By Lemma \ref{lemprim}, we may choose $(x,y)\in\bZ^2$ with 
$(x,y)\equiv (0,1)$ or $\equiv(1,0)\pmod{2}$ to $F_m(x,y)=d$ 
where $d$ is an odd divisor of $m^2+16$. 
Again by Lemma \ref{lemprim}, eight pairs $\pm(x,y)$, $\pm(y,-x)$ for $d$ and 
$\pm(x-y,x+y)$, $\pm(x+y,-x+y)$ for $-4d$ are all primitive solutions. 
Hence we may assume that $(x,y)\equiv(0,1)\pmod{2}$. 

These eight solutions correspond to the same $N$ as in Theorem \ref{th2} because 
\[
\frac{xy(x+y)(x-y)}{F_m(x,y)}
\]
is invariant under the actions of $x\longmapsto y\longmapsto -x$ and 
$x\longmapsto x-y$, $y\longmapsto x+y$. 

For $z=x/y$, the equation $(**)$ holds if and only if $f_A(z)=0$ 
where 
\[
A=\frac{mN+16}{-m+N}
\]
(see the proof of Theorem \ref{thC4}). 
We see that if $z=x/y$ is a root of $f_A(X)$ then other three roots of $f_A(X)$ are given by 
\[
\frac{z-1}{z+1}=\frac{x-y}{x+y},\quad 
-\frac{1}{z}=\frac{y}{-x},\quad 
-\frac{z+1}{z-1}=\frac{x+y}{-x+y}. 
\] 
Hence only the primitive solutions to $(*)$ which satisfy $(**)$ for $N$ 
are the eight solutions above. 
\end{proof}
\begin{corollary-nn}[Corollary \ref{corN}]
For $m\in\bZ\setminus\{0,\pm 3\}$, let $\mathcal{N}$ be the number of primitive solutions 
$(x,y)\in\bZ^2$ with $xy(x+y)(x-y)\neq 0$ to $F_m(x,y)=c$ where $c$ is a divisor of $4(m^2+16)$. 
Then we have 
\begin{align*}
\#\Bigl\{n\in\bZ\setminus\{\pm m\}\ \Big|\ L_n=L_m,\ n>0\Bigr\}=\frac{\mathcal{N}}{8}
\end{align*}
where $L_m=\mathrm{Spl}_\bQ f_m(X)$. 
In particular, if there does not exist $n\in\bZ\setminus\{\pm m\}$ with $L_n=L_m$ then 
$F_m(x,y)=c$ where $c$ is a divisor of $4(m^2+16)$ has only trivial solutions $(x,y)\in\bZ^2$ 
with $xy(x+y)(x-y)=0$. 
\end{corollary-nn}

By Theorem \ref{thC4}, we obtain 
\[
L_1=L_{103},\quad L_2=L_{22},\quad L_4=L_{956}. 
\]
Hence for $m\in\{1,2,4,22,103,956\}$, we get non-trivial eight primitive solutions 
to $(*)$ via Theorem \ref{th2} as on the following table: 
\begin{center}
{\rm Table} $2$\vspace*{5mm}\\
{\small 
{\renewcommand\arraystretch{1.2}
\begin{tabular}{|c|c|c|c|c|l|}\hline
$m$ & $N$ & $F_m(x,y)=c$  & $m^2+16$ & $xy(x+y)(x-y)$ & $(x,y)$ \\ \hline 
$1$ & $103$ & $-1$ & $17$ & $-6$ & $\pm(-2,1)$, $\pm(1,2)$\\ \hline
$1$ & $103$ & $4$ & $17$ & $24$ & $\pm(3,1)$, $\pm(-1,3)$\\ \hline
$2$ & $-22$ & $5$  & $20$ & $-6$ & $\pm(-2,1)$, $\pm(1,2)$\\ \hline
$2$ & $-22$ & $-20$ & $20$ & $24$ & $\pm(3,1)$, $\pm(-1,3)$\\ \hline 
$4$ & $-956$ & $1$ & $32$ & $-30$ & $\pm(-3,2)$, $\pm(2,3)$\\ \hline 
$4$ & $-956$ & $-4$ & $32$ & $120$ & $\pm(5,1)$, $\pm(-1,5)$\\ \hline 
$22$ & $-2$ & $125=5^3$ & $500=2^2 5^3$ & $-6$ & $\pm(-2,1)$, $\pm(1,2)$\\ \hline 
$22$ & $-2$ & $-500=-2^2 5^3$ & $500=2^2 5^3$ & $24$ & $\pm(3,1)$, $\pm(-1,3)$\\ \hline 
$103$ & $1$ & $-625=-5^4$ & $10625=5^4 17$ & $6$ &  $\pm(2,1)$, $\pm(-1,2)$\\ \hline 
$103$ & $1$ & $2500=2^2 5^4$ & $10625=5^4 17$ & $-24$ 
& $\pm(-3,1)$, $\pm(1,3)$\\ \hline 
$956$ & $-4$ & $28561=13^4$ & $913952=2^5 13^4$ & $-30$ 
& $\pm(-3,2)$, $\pm(2,3)$\\ \hline 
$956$ & $-4$ & $-114244=-2^2 13^4$ & $913952=2^5 13^4$ & $120$ 
& $\pm(5,1)$, $\pm(-1,5)$\\ \hline 
\end{tabular}
}}
\end{center}
\vspace*{3mm}
\section{Reducible case}\label{seRed}

In the reducible case $m\in\{0,\pm 3\}$, $f_m(X)$ splits as $f_0(X)=(X^2+2X-1)(X^2-2X-1)$ 
and $f_{\pm 3}(X)=(X^2\pm X-1)(X^2\mp 4X-1)$ over $\bQ$, and hence 
$\mathrm{Spl}_\bQ f_0(X)=\bQ(\sqrt{2})$ and $\mathrm{Spl}_\bQ f_{\pm 3}(X)=\bQ(\sqrt{5})$. 

If $m=0$, the trivial solutions correspond to $N=\pm m=0$. 

If $m=3$, then the eight trivial solutions $\pm(0,1)$, $\pm(1,0)$ for $c=1$ and 
$\pm(-1,1)$, $\pm(1,1)$ for $c=-4$ give $N=3$ and non-trivial eight solutions 
$\pm(2,1)$, $\pm(-1,2)$ for $c=-25$ and $\pm(-3,1)$, $\pm(1,3)$ for $c=100$ give $N=-3$. 

\begin{center}
{\rm Table} $3$\vspace*{5mm}\\
{\small
{\renewcommand\arraystretch{1.2}
\begin{tabular}{|c|c|c|c|c|l|}\hline
$m$ & $N$  & $F_m(x,y)=c$ & $m^2+16$ & $xy(x+y)(x-y)$ & $(x,y)$ \\ \hline 
$3$ & $-3$ & $-25$ & $25$ & $6$ & $\pm(2,1)$, $\pm(-1,2)$\\ \hline
$3$ & $-3$ & $100$ & $25$ & $-24$ & $\pm(-3,1)$, $\pm(1,3)$\\ \hline
\end{tabular}
}}
\end{center}
\section{Computational result}\label{seComp}

We do not know non-trivial primitive solutions to $(*)$ for $m\geq 0$ 
except on Table $2$ and Table $3$. 
By the correspondence as in Theorem \ref{th2}, in order to find primitive solutions 
to $(*)$ we should get $L_m=L_n$ for some $m\neq\pm n$. 
In \cite[Example 5.4]{HM09b}, however, we checked with the aid of computer that for integers 
$0\leq m<n\leq 10^5$, $L_m=L_n$ if and only if $(m,n)\in \{(1,103),(2,22),(4,956)\}$. 
Using Magma \cite{BCP97}, we can get the following: 
\begin{theorem}
For $0\leq m\leq 1000$, all non-trivial primitive solutions $(x,y)\in\bZ^2$, 
i.e. $xy(x+y)(x-y)$ $\neq 0$ and $\mathrm{gcd}(x,y)=1$, 
to $(*)$ are given on Table $2$ and Table $3$. 
In particular, for $0\leq m\leq 1000$ with $m\not\in\{1,2,4,22,103,956\}$ and $n\in\bZ$, 
$L_m=L_n$ implies $m=n$ or $m=-n$. 
\end{theorem}
%


\end{document}